\documentstyle[11pt,leqno,here]{article}
\input amssym.def
\input amssym.tex
\input psfig

\newcommand{\beql}[1]{\begin{equation}\label{#1}}
\newcommand{\eeq}{\end{equation}}

\makeatletter
\@addtoreset{figure}{section}
\def\thefigure{\thesection.\@arabic\c@figure}
\@addtoreset{table}{section}
\def\thetable{\thesection.\@arabic\c@table}

\def\@sect#1#2#3#4#5#6[#7]#8{\ifnum #2>\c@secnumdepth
     \def\@svsec{}\else
     \refstepcounter{#1}\edef\@svsec{\csname the#1\endcsname.\hskip .75em }\fi
     \@tempskipa #5\relax
      \ifdim \@tempskipa>\z@
        \begingroup #6\relax
          \@hangfrom{\hskip #3\relax\@svsec}{\interlinepenalty \@M #8\par}%
        \endgroup
       \csname #1mark\endcsname{#7}\addcontentsline
         {toc}{#1}{\ifnum #2>\c@secnumdepth \else
                      \protect\numberline{\csname the#1\endcsname}\fi
                    #7}\else
        \def\@svsechd{#6\hskip #3\@svsec #8\csname #1mark\endcsname
                      {#7}\addcontentsline
                           {toc}{#1}{\ifnum #2>\c@secnumdepth \else
                             \protect\numberline{\csname the#1\endcsname}\fi
                       #7}}\fi
     \@xsect{#5}}
\def\@begintheorem#1#2{\it \trivlist \item[\hskip \labelsep{\bf #1\ #2.}]}
\def\section{\@startsection {section}{1}{\z@}{-3.5ex plus -1ex minus
 -.2ex}{2.3ex plus .2ex}{\normalsize\bf}}

\makeatother

\begin{document}
\begin{center}
{\Large {\bf Improving Dense Packings of Equal Disks in a Square}}\\
\end{center}

\vspace{0.5\baselineskip}
\hspace{-.4in}
{\begin{tabular}{llll}
{\em David W. Boll} &
{\em Jerry Donovan} &
{\em Ronald L. Graham} &
{\em Boris D. Lubachevsky} \\
Hewlett-Packard &
Hewlett-Packard &
University of California &
Lucent Technologies\\
700 71st Ave &
700 71st Ave &
at San Diego &
Bell Laboratories \\
Greeley, CO 80634 &
Greeley, CO 80634 &
La Jolla, CA 92093 &
Murray Hill, NJ 07974 \\
david\_boll@hp.com &
jerry\_donovan@hp.com &
graham@ucsd.edu &
bdl@bell-labs.com
\end{tabular}

\setlength{\baselineskip}{0.995\baselineskip}

\begin{center}

{\bf ABSTRACT}
\end{center}

We describe a new numerical procedure for generating
dense packings of disks and spheres inside
various geometric shapes. We believe that in some of
the smaller cases, these packings are in fact optimal.
When applied to the previously studied cases of
packing $n$ equal disks in a square,
the procedure confirms all the
previous record packings
[NO1] [NO2] [GL],
except for
$n =$ 32, 37, 48, and 50 disks, where
better packings than those previously recorded are found.
For $n =$ 32 and 48,
the new packings are minor variations of the previous record
packings.
However, for $n =$ 37 and 50, the new patterns differ substantially.
For example, they are mirror-symmetric,
while the previous record packings
are not.


\setlength{\baselineskip}{1.2\baselineskip}
\section{Introduction}
\hspace*{\parindent}
We consider the task of arranging without overlaps a given number $n$
of congruent circular disks entirely inside 
a given square on the plane so that
the disks have the largest 
possible
diameter.
This disk packing problem has 
an equivalent formulation where one seeks to spread $n$ points 
(the centers of the disks) inside
a unit square so that the minimum point-to-point 
distance, usually denoted by $m = m_n$, is as large as possible.

In this paper, we describe a new experimental
approach and apply it to generate
new packings which are better than any previously known
[GL] [NO1] [NO2] [NO3]
for several values of $n$.
However, we do not prove the
optimality of these new packings,
although we suspect that some of them may indeed be optimal.
Clearly, a future goal would be to prove their optimality
(which becomes increasingly difficult as $n$ gets larger)
as is done in [NO3].

\section{New packings}
\hspace*{\parindent}
Figures ~\ref{t32-37} and ~\ref{t48-50} display 
the improved packings we found
and the corresponding old records for comparison.
For each number $n$ of disks presented,
the improvement in packing quality is small,
the first three digits of the displayed value of $m$ being identical.
Visually, the change in the packing pattern is also small
for $n=32$ and 48.
However, for $n=37$ and 50,
the difference in the corresponding packing patterns is substantial.
In particular, the new packings are mirror-symmetric,
while the old records are not.

\begin{figure}[H]
\centerline{\psfig{file=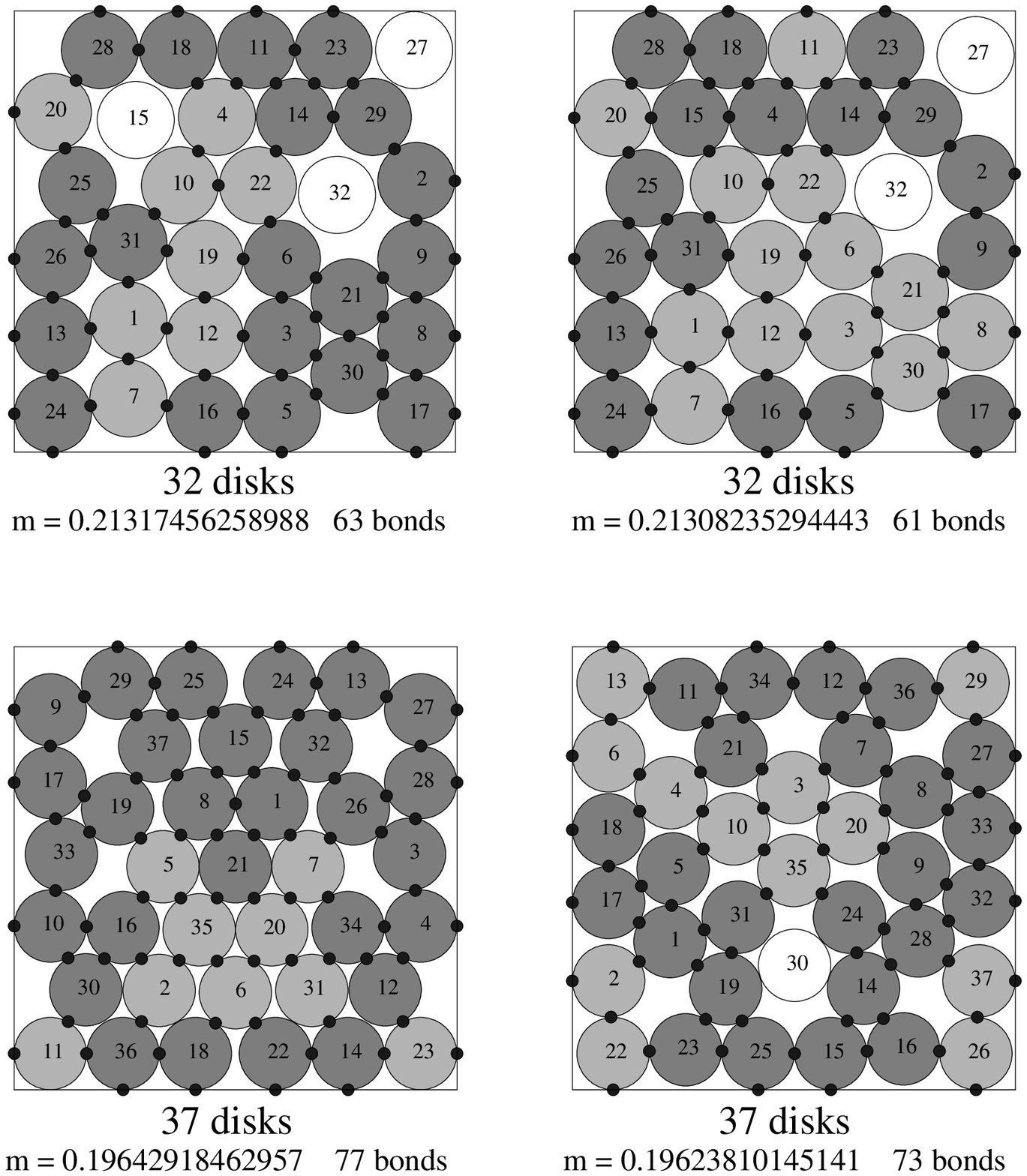,width=6.7in}}
\caption{{\it Left}. Improved packings of 32 and 37 disks in a square.
{\it Right}. Previous record packings of 32 and 37 disks in a square.
}
\label{t32-37}
\end{figure}

\begin{figure}[H]
\centerline{\psfig{file=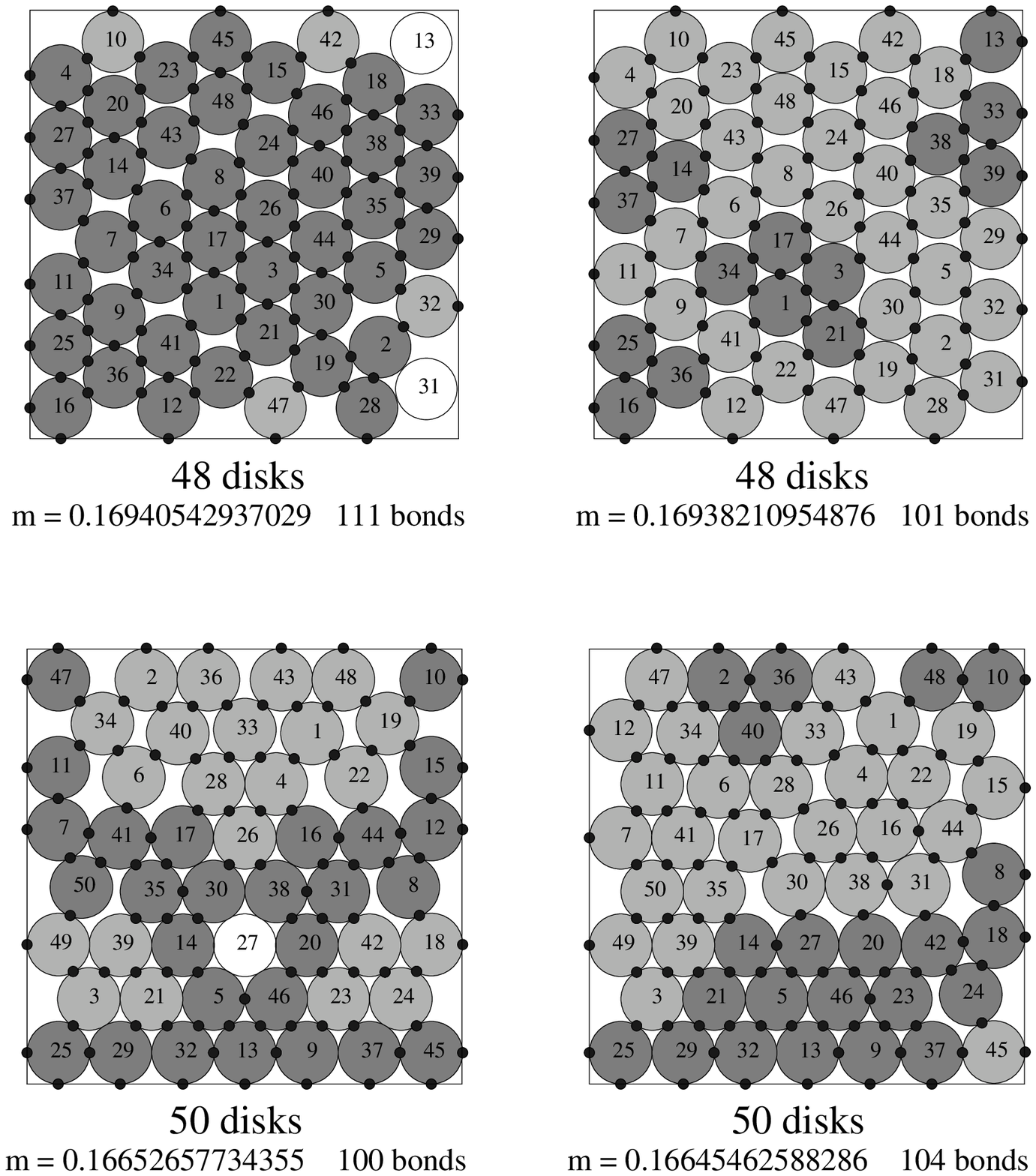,width=7in}}
\caption{{\it Left}. Improved packings of 48 and 50 disks in a square.
{\it Right}. Previous record packings of 48 and 50 disks in a square.
}
\label{t48-50}
\end{figure}

\clearpage

Each disk in a $n$-disk packing diagram
is provided with a unique label from 1 to $n$.
The labeling is arbitrary, although when possible we assign
the same label to the two disks
occupying similar positions in both the
improved and the previous record packings.

Black dots on the packing diagrams indicate
contacts, the so-called ``bonds.''
The total bond count is also provided with each packing.
A presence of a bond implies 
that the corresponding disk-disk or disk-boundary
distance is exactly zero, while the absence of a bond
in the spot of an apparent contact indicates a strictly positive distance.
As an example, here is a list of such ``almost'' contacts
in the improved packing of 37 disks 
shown at the bottom left in Fig.~\ref{t32-37}.
The distances are measured as fractions of the disk diameter.
\\\\
\begin{tabular} { r l r l r l }
21-5,7: &0.663349E-05 &20-35: &0.132670E-04 &22-31, 2-18: &0.110690E-03 \\
\\
bottom-11, 23:  &0.255689E-03 &2-30, 12-31: &0.128600E-01 & 30-36, 12-14:  &0.130805E-01
\end{tabular}
\\\\
Another example is the pair of disks 30 and 38
in the improved packing of 50 disks
shown at the bottom left in Fig.~\ref{t48-50}.
The distance between these disks is positive
(equal to 0.272038E-02 of the disk diameter),
and as a result disk 27 is a ``rattler'',
i.e., it is free to move within the ``cage'' formed
by its seven immobile neighbors, 
disks 38, 20, 46, 5, 14, 30, and 26.

All disks which are fixed in their positions are shaded,
the ``rattlers'' are unshaded.
We emphasize a packing structure by shading more heavily 
the disks that are ``more immobile'' than the others.
A disk $a$ belongs to such a more heavily shaded set,
if it is a member of a 3-clique ($a, b, c$), 
with each of $b$ and $c$ being either a disk or the bounding square side,
so that any two elements in the clique are in contact.
Thus, corner disk 45 in the improved 50-disk packing
(Fig.~\ref{t48-50}, bottom left)
is shaded more heavily because it has bonds with two sides of the square
which are incident to each other. In other words,
disk 45 belongs to 3-clique (disk 45, bottom side, right side).
In the same packing,
disk 15 belongs to 3-clique  (disk 12, disk 15, right side),
and disk 44 belongs to 3-clique (disk 8, disk 12, disk 44).
Notice that each improved packing 
increases the number of such ``more immobile'' disks 
compared to the number in the previous record,
a phenomenon which is intuitively reasonable,
but must not necessarily hold in general.

\section{How does one know a packing exists?}
\hspace*{\parindent}
We hope that the improved packings we present are optimal,
i.e., the densest
possible for their value of $n$.
However, 
the only evidence we have for this
is that no one has been able (so far)
to find better packings.
The only values of $n$ for which optimality has
been proved are $n \le 27$, and these proofs involve substantial
computation as the value of $n$ increases,
and the examination of 
an exponentially increasing number of cases.

Before proving a displayed packing is optimal, 
it is first necessary to prove that it actually exists!
That is,
that there is a geometrical configuration corresponding
to the points and edges of a given diagram in which all 
inter-center
distances are at least one disk diameter
and center-boundary distances are at least a disk radius.
For example,
existence of the packings in Figures ~\ref{t32-37} and ~\ref{t48-50},
is not obvious (to us!) from just looking
at the diagrams.
(For simplicity of discussion,
we will assume for now that all disk diameters are 1,
and we are trying to minimize the side-length of the smallest
square into which $n$ nonoverlapping disks
can be packed.)

While we do not offer such a proof of existence for our
four new packings,
the following procedure could (in principle)
be used to generate one.
With
$z_1 ,z_2 ,...z_n$ 
denoting the centers of the disks
in a 
given a packing diagram (with $z_i = ( x_i ,y_i )$),
we can generate a (usually overdetermined)
set of equations in the $x_i$ and $y_i$ expressing
the various bonds (or edges) in the diagram.
We can then eliminate all the variables except one
(say $x_1$) from the equations (using resultants),
to obtain a single polynomial
equation $P(x_1) = 0$ which must be satisfied by any solution
$(x_1, y_1,...,x_n , y_n)$ of the original set of equations.
We can then express each $z_i$ in terms of $x_1$,
and finally bound the distances between every pair $z_i$ and $z_j$
and between every $z_i$ and the boundary.
Those which have bonds correspond to edges which have exactly 1
or where the distance of a disk to the boundary is exactly 0.
Those which do not must have inter-center distance strictly greater
than 1 or the disk-boundary distance must be strictly positive,
and this will be confirmed
by using sufficient accuracy in approximating
the appropriate root of $P(x_1) = 0$.
As is typical in this type of variable elimination,
care must be taken in deciding just how it is to be performed,
since there is a strong tendency for the size and number of terms
to grow beyond control if
one is not careful!

A more practical, but, perhaps, less convincing variant of
the above procedure 
starts as above with the set of equations with respect to $z_i$.
Instead of eliminating all the variables, but one,
the set is transformed into a quadratic minimization
problem and a gradient iterative procedure is employed
to find the minimum.
The packing diagram itself may be used to input an initial approximation
for the solution. 
The iterative calculations can be 
implemented with arbitrary precision.
If the procedure numerically converges with high precision,
say using 100 decimal digits,
yielding the minimum in which
for each disk-disk bond the inter-center distance is 1
and for each disk-boundary bond the distance is 1/2,
say using 100 decimal digits, it would be an evidence
of the packing's existence.

Note that if we just had a diagram 
in which only the disks (or disks centers)
were shown, 
without the proposed bonds (edges),
then a large (potentially exponential with $n$)
number of guesses would have to be made as to where
bonds might be,
before applying either of the aforementioned
procedures to each guess.

Another method (also not a proof!)
of convincing oneself of a
packing's existence, its {\em rigidity},
and of the correctness of the presented parameters
is to apply the procedure described in detail
in the next section.
The procedure consists of two Phases.
The novel part is Phase 1. It delivers
an approximate configuration.
In difficult packing cases, like those in
Figures ~\ref{t32-37} and ~\ref{t48-50}
the degree of approximation does not suffice
to identify bonds and to determine a sufficient number of decimal digits
of packing parameters, e.g., the value of $m$.
Phase 2 takes the configuration generated
in Phase 1 as a starting point
and runs the ``billiards'' simulation procedure
described previously (see, e.g, [L]).
The latter usually provides a clear-cut bond identification,
e.g., all pairs of disks with pairwise distances less than $10^{-11}$
of the disk diameter form bonds, while the other pairwise
distances are larger than $10^{-5}$ of the disk diameter.

For the sake of convincing
ourselves of the packing's existence
and its rigidity,
we 
perturbed
the path to the final configuration.
For example, we may turn the intermediate
configuration after Phase 1 by $90^{\circ}$.
After applying Phase 2, the resulted configuration
should be (and was) the $90^{\circ}$ rotation of
the packing produced in the unperturbed path.
The
identified
bonds as well as the resulting packing parameters
(with high precision)
were
independent of the perturbation.

\section{Why a new procedure?}
\hspace*{\parindent}
Note that the ``billiards'' algorithm alone
is also able to converge to the record packing
after many attempts with different
randomly chosen starting points.
It seems that, at least for packing
disks in a square, the combination of the two
procedures
performs substantially better.
For example, more than 1000 attempts were needed 
by the ``billiards'' algorithm to generate a single
record packing of 32 disks in Fig. \ref{t32-37}.
For 37 disks more than 5000 attempts were needed.
The new algorithm 
converged in under 30 attempts for 32 disks
and under 100 attempts for 37 disks.

We observed that
each next digit of precision takes roughly several times
more computing time
than the previous digit.
Thus, computing $k$ digits takes time
increasing roughly exponentially with $k$.
This observation holds, of course, only as long
as the procedure {\em is able} to deliver new digits.
When all computations are performed with double precision
(i.e., with error less than $10^{-14}$),
the ``billiards'' algorithm usually delivers more than 13 digits,
while the new Phase 1 algorithm becomes stuck and
terminates at about the 10th digit.
Thus,
one run of Phase 1 terminates much sooner
(yielding an approximate packing)
than one run of the
``billiards'' algorithm 
(which yields an exact packing).
Because of that,
the difference between the two
procedures goes from months
of computing to just a few hours.
A substantial share of the latter's time
is spent in the single finalizing
run of Phase 2 (the ``billiards'' algorithm), where
the starting point is 
the record packing obtained
by multiple runs of Phase 1.

Conjectured optimal packings of $n$ disks 
in a square 
for various $n \neq 32,37,48$ or 50 (see [GL])
have also been 
found by the new combination algorithm,
which also ran more quickly in these cases
than the ``billiards'' algorithm alone.

\section{The new packing algorithm}
\hspace*{\parindent}
We describe the algorithm of Phase 1
in a form applicable for packing equal spheres
in variously-shaped bounded regions of euclidean $n$-space.
Our computer realization of this procedure packs
equal disks in rectangles, equilateral triangles, or circles.
We have used it primarily for packing disks in squares.

\clearpage

\begin{enumerate}

\item Arbitrarily initialize the center position of each sphere
so that the spheres do not intersect while each has a unit radius.
Arbitrarily initialize a direction of motion for each sphere.
(In our packing-disks-in-a-square experiments, we usually randomly scatter
the disks in a sufficiently large square and for a disk centered at $a$,
we set its initial motion direction to be $-a$.)

\item Determine the smallest size bounding shape 
that contains all spheres. 
(Specifically, we determine the smallest bounding square for all disks.
Without losing packing performance,
the square size minimization is done more easily if restricted
to the squares centered at the point (0,0)
with sides parallel to the x and y-axes.)

\item Set the move step $s$ to its initial predefined value, $s = s_0 > 0$.
(We usually take $s_0 = 0.25$.)

\item Initialize the counter $impatience = 0$.

\item Execute the following Steps 6 - 8 repeatedly
while the move step $s$ remains larger than a predetermined
small threshold $\epsilon$. Once $s \le \epsilon$, terminate the run.
(We usually take $\epsilon = 10^{-10}$.)

\item Select the spheres sequentially without repetition in some order
and attempt to move each sphere a distance $s$ from its current position
using procedure $move\_attempt$ detailed below. Each $move\_attempt$
for a sphere returns either $success$ or $failure$.
(In our experiments, the order of sphere selection appears 
non-essential. We set the  selection order randomly at the start of the procedure
and randomly reset it each time when resetting $s$ in Step 8 below.)

\item Determine the smallest size bounding shape 
that contains all spheres.
If the new size of the shape is not smaller than the previous size,
then increment $impatience$ by 1.

\item If $impatience$ exceeds a predefined threshold $tolerance$
or if all recent $move\_attempt$s made in Step 6 
(their number is equal to the number of spheres) fail, then
reset the move step to a smaller value, $s \leftarrow s * q$,
where $q < 1$ is a predefined constant. If $s$ is thus changed,
then reset $impatience \leftarrow 0$.
(While packing disks in a square, 
we usually take $tolerance = 1000$ and $q = 0.43$.)

\end{enumerate}

    Procedure $move\_attempt$.

\begin{itemize}

\item Attempt to move the sphere a distance $s$ 
along its previously chosen direction.
(The initial direction has been chosen in Step 1 of the main procedure.)

\item If in the new position the sphere does not intersect other spheres
and does not extend beyond the region boundary, then $move\_attempt$ terminates after
accepting the new sphere position and reporting $success$ to the calling
procedure.

\item Otherwise start with the old sphere position again and
reset the direction of motion using
the procedure  $reset\_direction$ detailed below and attempt to move
the sphere a distance $s$ in the new chosen direction.

\item If in the new position the sphere does not intersect other spheres
and does not extend beyond the region boundary, then $move\_attempt$ terminates after
accepting the new sphere position and reporting $success$ to the calling
procedure.

\item Otherwise $move\_attempt$ terminates after
reporting $failure$ to the calling procedure.
The old sphere position does not change.

\end{itemize}

    Procedure $reset\_direction$.

\begin{itemize}

\item The directional vector is formed as the
sum of the ``repelling'' vectors, summed over all the obstacles that can be reached
by the given sphere in a single straight motion of length $s$.

\item Another sphere yields 
a repelling obstacle for the given sphere if the center 
of the former sphere is located
at a distance $d_1 \le 2(1+s)$ from the given sphere center. 
The obstacle is the point on the former sphere 
closest to the given sphere.

\item A point on the region boundary qualifies as a repelling obstacle
if it yields
a local minimum of the distance 
from the center of the given sphere to the boundary
and if
the distance
from the point to the given sphere center is $d_2 \le 1+s$.
(Thus, the region boundary may yield multiple repelling obstacles.
For example, 
the two sides that form a square corner
yield two boundary obstacles 
for a disk located nearby the corner.)

\item The direction for a repelling vector is formed 
by connecting the obstacle
point with the center of the given sphere toward the center.
(Thus, another sphere $b$ repels the given sphere $a$
in the direction from the center of $b$ to the center of $a$.
A boundary obstacle repels in the direction from
the point of local minimum of the distance to the center of $a$.)

\item The length of a repelling vector is $d^\alpha$,
where $\alpha$ is a preselected parameter, and
$d = d_1$ if the repelling obstacle is given by a sphere,
and $d = 2 d_2$ if the repelling obstacle is given by the boundary.
(Most of our disk packing experiments were done
with $\alpha = 1$.)

\item $reset\_direction$ returns this
computed directional vector if it is non-zero.
Otherwise, it returns a default directional vector.
(The choice of the default directional vector is not
critical for the packing performance
as long as the vector does not usually degenerate to 0.
When packing disks in a square, we choose $-a$ as the default
directional vector for the given circle centered at $a$.)

\end{itemize}

\section{Concluding remarks}
\hspace*{\parindent}
We have found that combining two heuristic procedures for producing
packings gives a substantial improvement over the use of either
one alone,
a phenomenon which often occurs in approximation algorithms.
It may well be that additional heuristics can improve the performance
even further (and we are currently experimenting with such algorithms).
The next challenge of course will be to prove that the packings produced 
are in fact optimal
(and also unique, which we also believe, up to the position
of the so-called ``rattlers'').
A good start in this direction has been made in [NO3].
Even more ambitious would then be to prove optimality
for infinite families of packings, as can be done when packing 
$n(n+1)/2$ equal disks in an equilateral triangle.


\begin{thebibliography}{MMMM}
\bibitem[GL]{GL}
R.~.L. Graham and B.~D. Lubachevsky,
Repeated Patterns of Dense Packings of Equal Disks in a Square,
{\em Electronic Journal of Combinatorics} {\bf 3(1)} (1996), \#R16.
\bibitem[L]{L}
B.~D. Lubachevsky,
How to simulate billiards and similar systems, {\em J. Computational
Physics} {\bf 94} (1991), 255--283.
\bibitem[NO1]{NO1}
K.~J. Nurmela and P.~R.~J. \"{O}sterg{\aa}rd,
Packing up to 50 equal circles in a square,
{\em Discrete Computational Geom.} {\bf 18} (1997), 111-120.
\bibitem[NO2]{NO2}
K.~J. Nurmela and P.~R.~J. \"{O}sterg{\aa}rd,
Optimal Packings of Equal Circles in na Square,
in: {\em Combinatorics, Graph Theory, and Algorithms},
Vol. II, Y. Alavi, D.R. Lick, and Schwenk (eds.),
New Issues Press, Kalamazoo 1999, pp.671-680.
\bibitem[NO3]{NO3}
K.~J. Nurmela and P.~R.~J. \"{O}sterg{\aa}rd,
More optimal packings of equal circles in a square,
{\em Discrete Computational Geom.} {\bf 22} (1999), 439-457.
\end{thebibliography}
\end{document}